\documentclass[letterpaper,10 pt,conference]{ieeeconf}
\IEEEoverridecommandlockouts
\usepackage{amsmath,amssymb,amsfonts}
\usepackage{algorithmic}
\usepackage[ruled]{algorithm2e}
\usepackage{graphicx}
\usepackage{textcomp}
\usepackage{xcolor}
\usepackage{cite}
\def\BibTeX{{\rm B\kern-.05em{\sc i\kern-.025em b}\kern-.08em
    T\kern-.1667em\lower.7ex\hbox{E}\kern-.125emX}}
    
\begin{document}

\title{Learning-Accelerated ADMM for Distributed DC Optimal Power Flow
\thanks{This work was authored by the National Renewable Energy Laboratory, operated by Alliance for Sustainable Energy, LLC, for the U.S. Department of Energy (DOE) under Contract No. DE-AC36-08GO28308. Funding was provided by the Autonomous Energy Systems project funded by the National Renewable Energy Laboratory's Laboratory Directed Research and Development program. The views expressed in the article do not necessarily represent the views of the DOE or the U.S. Government. The U.S. Government retains and the publisher, by accepting the article for publication, acknowledges that the U.S. Government retains a nonexclusive, paid-up, irrevocable, worldwide license to publish or reproduce the published form of this work, or allow others to do so, for U.S. Government purposes.}
}

\author{
    \authorblockN{
        David Biagioni\authorrefmark{1},
        Peter Graf\authorrefmark{1},
        Xiangyu Zhang\authorrefmark{1},
        Ahmed Zamzam\authorrefmark{1},
        Kyri Baker\authorrefmark{2},
        Jennifer King\authorrefmark{1}}
    \authorblockA{\\
        \authorrefmark{1}National Renewable Energy Laboratory, 15013 Denver West Parkway, Golden, CO, 80401\\
        Email: \{dave.biagioni, peter.graf, xiangyu.zhang, ahmed.zamzam, jennifer.king\}@nrel.gov}
    \authorblockA{\\
        \authorrefmark{2}Dept. of Civil, Environmental and Architectural Engineering\\
        University of Colorado at Boulder, UCB 428, Boulder, CO, 80309\\
        Email:  Kyri.Baker@colorado.edu}
}

\maketitle

\begin{abstract}
We propose a novel data-driven method to accelerate the convergence of Alternating Direction Method of Multipliers (ADMM) for solving distributed DC optimal power flow (DC-OPF) where lines are shared between independent network partitions.  Using previous observations of ADMM trajectories for a given system under varying load, the method trains a recurrent neural network (RNN) to predict the converged values of dual and consensus variables.  Given a new realization of system load, a small number of initial ADMM iterations is taken as input to infer the converged values and directly inject them into the iteration.  We empirically demonstrate that the online injection of these values into the ADMM iteration accelerates convergence by a significant factor for partitioned 14-, 118- and 2848-bus test systems under differing load scenarios. The proposed method has several advantages: it maintains the security of private decision variables inherent in consensus ADMM; inference is fast and so may be used in online settings;  RNN-generated predictions can dramatically improve time to convergence but, by construction, can never result in infeasible ADMM subproblems; it can be easily integrated into existing software implementations. While we focus on the ADMM formulation of distributed DC-OPF in this paper, the ideas presented are naturally extended to other distributed optimization problems.
\end{abstract}


\section{Introduction}\label{sec-intro}
The electric power grid is continually progressing towards a more complex, uncertain, and decentralized state.  This fact stems from a variety of sources including higher penetration of renewable generation, increased presence of smart devices and subsystems, and market deregulation, to name a few.  While this progression presents a number of operational and analytical challenges, the corresponding increase in available data paves the way for new approaches to improving system-wide efficiency and coordination.

In this paper, we consider the specific operational problem of direct current optimal power flow (DC-OPF) in which the network is decomposed into independently operating partitions.  Because each partition is connected to others via one or more branches, agreement on these line flows is required as part of an optimal solution.  While there are many possible ways to solve this \emph{consensus} problem in a distributed fashion, we focus on the Alternating Direction Method of Multipliers (ADMM)  due to its recent popularity for solving such problems \cite{bertsekas1989parallel,boyd2011distributed,chang2014multi,zhang2014asynchronous}.  We note, however, that the ideas proposed in this paper are readily extensible to other iterative solution techniques.

In the spirit of several recent analyses of ADMM \cite{nishihara2015general,francca2018admm}, we view the iterations as the time steps of a discrete, stable dynamical system whose equilibria correspond to an optimal solution of the underlying optimization problem.  The proposed method, which we call \emph{learning-accelerated ADMM (LA-ADMM)}, aims to leverage previously observed trajectories as training data to build a machine learning (ML) model to predict the equilibrium values using only a small number of initial iterates.  Once such a model is trained, optimal values can be inferred online and injected directly into the iteration to accelerate convergence to a feasible and optimal solution.  During prediction, we use a recurrent neural network (RNN) based on the Gated Recurrent Unit (GRU) \cite{cho2014learning} and argue that this is a sound choice based on the connection between ADMM and numerical schemes for solving a particular dynamical system  \cite{nishihara2015general,francca2018admm}.

The use of ML in power systems is a growing area of research. Many straightforward applications of ML in OPF and related problems are either \emph{centralized} in the sense that training requires complete knowledge of the system's state and dynamics \cite{Zamzam_learn_19,chatzos2020highfidelity,zamzam2019data,yang2020robust}, or \emph{intrusive} in the sense that target models are created \emph{de novo} or to replace an existing one; see, e.g., \cite{karagiannopoulos2019data,dobbe2019towards,xavier2019learning,yang2019two,zamzam2019energy} for recent examples. Our approach requires neither centralization nor intrusiveness. In contrast to many works using ML for OPF, we use ML not to obtain the optimal solutions of the OPF itself; rather, we use ML to accelerate the process of solving the OPF. Importantly, this avoids infeasibility of the final solution, unlike ML+OPF methods which have no feasibility guarantees or require a feasibility post-processing step to ensure power flow constraints hold \cite{Zamzam_learn_19, DeepOPF, zhao2020deepopf+}. While the proposed method can leverage centralized models to accelerate the training process, it only requires access to the consensus variables for training and inference. Because predicted values are fed into the optimization models as iteration parameters, the method is distinctly non-intrusive.

In addition to the interest in ML increasing in the power systems community, interest in distributed optimization algorithms is also on the rise. This is partially due to the increasing prevalence of physically distributed, autonomous systems and the frequent intractability of centralized formulations that assume complete knowledge of the system's state and dynamics. The optimal power flow problem, in particular, is amenable to the use of distributed methods as evidenced by several recent reviews \cite{kargarian2016toward,molzahn2017survey} and applications utilizing ADMM  \cite{erseghe2014distributed,peng2014distributed,scott2015distributed,wang2016fully,ma2016consensus,zhang2017dynamic,zhang2017distributed}. This paper combines ML with ADMM to leverage both the distributed nature of ADMM and the decrease in convergence time from using ML to estimate the dual and consensus variables.

The paper is organized as follows.  In section \ref{sec-methods}, we briefly summarize the components of the proposed method, namely, the DC-OPF problem, its distributed formulation via ADMM, and the RNN used for the prediction task.  In section \ref{sec-experiments}  we describe empirical experiments for three test systems and present the results of using standard ADMM versus LA-ADMM.  Section \ref{sec-discussion} discusses 
further issues, conclusions and future work.

\begin{figure}[t!]
\centerline{\includegraphics[width=0.95\linewidth]{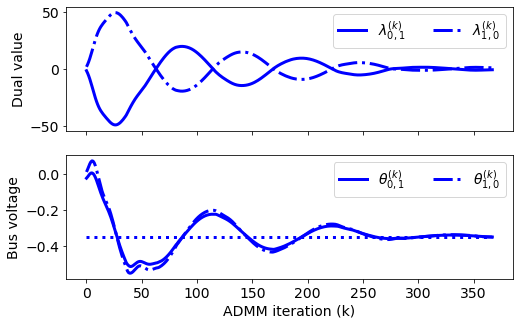}}
\caption{Convergence of the DC-OPF ADMM iteration for the voltage value at a boundary bus in the 14-bus test system.  The top panel shows the dual variables for this bus for each partition.  The bottom panel shows the value of the shared voltage for the two partitions, with the optimal consensus value as a thin horizontal line.  This figure illustrates the viewpoint of ADMM as a discrete dynamical system approaching a steady state at optimality as discussed in Section \ref{subsec-rnn}.}
\label{fig-admm-iteration}
\end{figure}

\section{Methodology}\label{sec-methods}

\subsection{DC Optimal Power Flow Formulation}

\newcommand{\gen}[0]{\mathbf{g}}
\newcommand{\ang}[0]{\boldsymbol{\theta}}

The aim of the DC-OPF is to determine the least-cost dispatch of generation that satisfies the load demand subject to power flow limits on the network lines. Using a DC approximation for the AC power flow equations leads to the following optimization problem:
\begin{align}
\operatornamewithlimits{min}_{\gen,\ang}\quad& f(\gen) \label{eqn-dcopf-objective}\\
\operatorname{s.t.}\quad& \mathbf{H} \ang + \gen - \mathbf{d} = \mathbf{0}, \label{eqn-power-balance}\\
&\underline{\gen} \leq \gen \leq \overline{\gen}, \qquad \underline{\mathbf{f}} \leq \mathbf{K} \ang \leq \overline{\mathbf{f}}, \qquad \theta_0 =0.
\label{eqn-flow-limits}
\end{align}
Here $\gen,\mathbf{d}\in\mathbb{R}^n$ are the active power generation and load at each of $n$ buses, respectively, and $\ang\in\mathbb{R}^n$ collects the phase angles at each of $n$ buses, with $\theta_0$ denoting the angle at the reference bus. The matrix $\mathbf{K} :=  \mathbf{B} \mathbf{A}$, where $A\in\mathbb{R}^{m\times n}$ is a network directed graph incidence matrix, where $m$ is the number of lines in the network, and the diagonal matrix $\mathbf{B} \in \mathbb{R}^{m\times m}$ collects the line susceptances on the diagonal. Also, the matrix $\mathbf{H} \in \mathbb{R}^{n \times n}$ is defined as $\mathbf{H} := \mathbf{A}^T \mathbf{B} \mathbf{A}$.
Upper and lower bounds on $\gen$ are denoted by $\overline{\gen}, \underline{\gen}$, respectively. Similarly, upper and lower line flow limits are denoted by $\overline{\mathbf{f}}$ and $\underline{\mathbf{f}}$, respectively. The vector of all zeros, $\mathbf{0}$, has dimension determined by context (here, $\mathbf{0}\in\mathbb{R}^n$). The cost function $f:\mathbb{R}^n \rightarrow \mathbb{R}$ representing the cost of generation is assumed to be linear, i.e.,
$f(\gen) = \mathbf{c}^T \gen$ with non-negative cost vector,  $\mathbf{c}\geq\mathbf{0}$.  

\begin{algorithm}[t!]
\SetKwInput{KwInput}{Input}                
\SetKwInput{KwOutput}{Output}              
\DontPrintSemicolon
\SetAlgoLined
\KwInput{DC-OPF ADMM problem specification; $K$-step RNN to predict converged $\boldsymbol{\lambda}_{su}^*,\overline{\ang}_{su}^*$}
\KwOutput{Distributed solution to DC-OPF}
 Set $k\leftarrow0$\;
 \While{ADMM not converged}{
   \uIf{$k=K$}{
      Predict converged values:  $\boldsymbol{\lambda}_{su}^*,\overline{\ang}_{su}^*$\;
      Overwrite dual variables: $\boldsymbol{\lambda}_{su}^{(k)} \leftarrow \boldsymbol{\lambda}_{su}^*$\;
      Overwrite consensus variables: $\overline{\ang}_{su}^{(k)}\leftarrow\overline{\ang}_{su}^*$\;
   }
   Execute ADMM iteration\;
   $k\leftarrow k+1$
 }
 \label{alg-la-admm}
 \caption{Learning Accelerated-ADMM}
\end{algorithm}

\subsection{Distributed DC-OPF Formulation via ADMM}\label{subsec-admm-dcopf}
In the distributed DC-OPF problem, the network is partitioned physically (i.e., into balancing areas) and/or computationally into subproblems that may share some number of the optimization variables.  In this paper, we specifically consider the partitioning of buses into $S$ disjoint sets indexed by $\mathcal{I}_s$, such that $\cup_{s=1}^S\mathcal{I}_s = \{1,...,n\}$.  In other words, a bus $i$ belongs to partition $s$ if and only if $i \in \mathcal{I}_s$. Such partitioning naturally leads to bus classification into buses that are connected only to buses inside the same partition, and buses that are connected to buses in other partitions. We denote the set of buses inside partition $u$ that are connected to buses in partition $s$ by $\mathcal{J}_{su}$, where $s, u$ refer to adjacent partitions. We use the terminology of \emph{public} and \emph{private} to distinguish between decision variables that are shared between partitions and those that are internal to a partition, respectively. Using the above notation, we note that public variables consist of phase angles at buses that are connected to buses in other partitions. Let the number of buses that are connected to buses outside their partitions be $n_{\text{pub}}$. In contrast, all generation decision variables and phase angles of buses that are only connected to buses inside the same partition are private decision variables. Thus, the number of these private variables is $n_{\text{pri}} = 2n - n_{\text{pub}}$, where $2n$ corresponds to both bus voltage and generator setpoints at each private bus. Finally, we use subscript $s$ to denote the partition membership of decisions, and use $\ang_{su}$ to denote the phase angles of buses in partition $u$ that are connected to buses in partition $s$.

The global DC-OPF problem \eqref{eqn-dcopf-objective}-\eqref{eqn-flow-limits} can be expressed in terms of the partitions by
\begin{align}
\operatornamewithlimits{min}_{\gen, \ang}\quad& \sum_{s=1}^S f_{s}(\gen_s) \label{eqn-dcopf-objective-dist}\\
\operatorname{s.t.}\quad& \mathbf{H}_s \ang_s + \sum_{u} \mathbf{H}_{su} \ang_{su} + \gen_{s} - \mathbf{d}_{s} = \mathbf{0}_{s},&& \forall{s} \label{eqn-power-balance-dist}\\
&\underline{\gen}_{s} \leq \gen_{s} \leq \overline{\gen}_{s}, \qquad &&\forall s\\
&\underline{\mathbf{f}}_s \leq \mathbf{K}_s \ang_s + \sum_u \mathbf{K}_{su} \ang_{su} \leq \overline{\mathbf{f}}_s, &&\forall{s} \label{eqn-flow-limits-dist}\\
& \boldsymbol{\theta}_{su} = {\bf E}_{su} \boldsymbol{\theta}_u &&\forall s, u \label{eqn-consensus}
\end{align}
where $\mathbf{H}_{su}$ denotes the sub-block of the matrix $\mathbf{H}$ that selects the rows  $\mathcal{I}_s$ and columns $\mathcal{J}_{su}$. In other words, the constraint enforces power balance on all internal nodes subject to both internal and shared decision variables.  The cost function is also partitioned over the network partitions where $f_{s}$ is function of only the generation at generation inside the $s$ partition. The matrix ${\bf K}$ is also partitioned such that for each partition the flow limits are respected for all lines connected to at least one bus in the partition. The matrices ${\bf E}_{su}$ is a selection matrix that selects the phase angles of buses in $\mathcal{J}_{su}$.

Were it not for the consensus constraints \eqref{eqn-consensus},
the above problem \eqref{eqn-dcopf-objective-dist}-\eqref{eqn-flow-limits-dist} could be trivially decomposed into $S$ disjoint problems whose independent solution yields the global minimum.  In the presence of these constraints, however, special handling is required to enforce these constraints.  ADMM accomplishes this task via an augmented Lagrangian formulation (see the classic reference \cite{boyd2011distributed} for more details) that drives local copies of public variables into consensus with the partitions that share them.  Algorithmically, consensus ADMM equates to the following iterative scheme consisting of $S$ independent primal optimizations followed by a centralized dual update,
\begin{align}
&\{{\bf g}_s^{(k+1)}, \boldsymbol{\theta}_s^{(k+1)}, \boldsymbol{\theta}_{su}^{(k+1)} \}= \operatornamewithlimits{argmin}_{{\bf g}_s, \boldsymbol{\theta}_s, \boldsymbol{\theta}_{su} \in \mathcal{C}_s} f_s(\gen_s)   \nonumber \\&\quad + \sum_{u} {\boldsymbol{\lambda}_{su}^{(k)}}^T \left( \ang_{su} - \overline{\ang}_{su}^{(k)} \right)
 +\frac{\rho}{2} \big|\big| \ang_{su} - \overline{\ang}_{su}^{(k)} \big|\big|_2^2
\label{eq-admm-xpdate}\\
&\overline{\ang}_{su}^{(k+1)} = \frac{1}{2} \left( \ang_{su}^{(k+1)} + {\bf E}_{su} \ang_u^{(k+1)} \right)
\label{eq-admm-mean-update}\\
&\boldsymbol{\lambda}_{su}^{(k+1)} = \boldsymbol{\lambda}_{su}^{(k)} + \rho \left( \ang_{su}^{(k+1)} - \overline{\ang}_{su}^{(k+1)} \right). \label{eq-admm-lambda-update}
\end{align}
where we have introduced Lagrange multipliers  $\boldsymbol{\lambda}_{su}\in\mathbb{R}^{|\mathcal{J}_{su}|}$ and used $\mathcal{C}_s$ to denote the local constraint set for partition $s$, i.e., constraints \eqref{eqn-power-balance-dist}-\eqref{eqn-flow-limits-dist}.

\subsection{Recurrent Neural Networks}\label{subsec-rnn}
As alluded to in Section \ref{sec-intro}, one way to view the ADMM iteration is as a numerical integration scheme for solving the dynamical system, sometimes referred to as \emph{gradient flow},
\begin{align}
\frac{d}{dt}u(t) = -\nabla f(u(t)).
\end{align}
Note that steady states of this system coincide with local optima of $f$ since $du/dt=0 \iff \nabla f = 0$.  This point of view has proven useful in several recent works analyzing the convergence and designing accelerated variants of ADMM \cite{nishihara2015general,francca2018admm}.  It is also empirically intuitive; see, e.g., Fig. \ref{fig-admm-iteration} showing convergence of a subset of ADMM variables for the 14-bus system described in Section \ref{sec-experiments}.

\begin{table}[t!]
\caption{Test systems. Notation: $n_{pub}$ denotes the number of boundary buses connecting two partitions; $|x_{su}|$ denotes the number of corresponding optimization variables (dual and primal) and $\rho$ is the ADMM step size.}
\begin{center}
\begin{tabular}{|l|c|c|c|c|c|c|}
\hline
\textbf{Network} & Buses & Branches & Partitions & $n_{pub}$ & $|x_{su}|$ & $\rho$ \\
\hline
IEEE 14 & 14 & 20 & 2 & 6 & 24 & 1\\
\hline
IEEE 118 & 118& 186 & 4 & 22 & 88 & 100\\
\hline
RTE 2848 & 2848 & 3776 & 18 & 257 & 1028 & 100\\
\hline
\end{tabular}
\label{tab-system-data}
\end{center}
\end{table}

In selecting a machine learning algorithm, we adopted the dynamical system interpretation of ADMM and only considered algorithms suitable for predicting sequential data.  While many choices exist for this task, we selected the GRU \cite{cho2014learning} because it is a well established architecture in deep learning, is relatively lightweight in terms of number of parameters and, as a result, is extremely easy implement and train using modern open source software.  While not formally reported in this paper, preliminary experiments also suggested that an RNN architecture was superior for this prediction task when compared with other off-the-shelf ML algorithms.

\subsection{LA-ADMM Algorithm}

The LA-ADMM algorithm is a modification of standard ADMM in which the iteration is interrupted once to inject predicted values for the converged dual and consensus variables.  We assume that an RNN has been previously trained using $K$ steps of ADMM for which the optimal solution is known, either by gathering data from fully convergent ADMM iterations or, where applicable, by solving a centralized version of the problem.  Thus, the RNN takes as input the values $\boldsymbol{\lambda}_{su}^{(k)},\overline{\ang}_{su}^{(k)}$ with $k=1,\dots,K$ and generates predictions for the converged values which we denote by  $\boldsymbol{\lambda}_{su}^*,\overline{\ang}_{su}^*$.  These predicted values overwrite the current iteration variables and become the \emph{de facto} parameters for partition subproblems.  Pseudocode for LA-ADMM is given in Algorithm \ref{alg-la-admm}.

A natural question is whether it would be beneficial to invoke the RNN prediction \emph{every} $K$ steps rather than just once.  While this seems like a reasonable idea, our experiments suggest that such a method suffers from accumulating prediction errors since, for all but the first prediction step, inputs to the RNN have been perturbed by previous prediction steps.  Empirically we observed that this led to poor outcomes for the algorithmic settings we considered and so did not investigate it further.

\section{Experiments and Results}\label{sec-experiments}

\begin{figure}[t!]
\centerline{\includegraphics[width=0.95\linewidth]{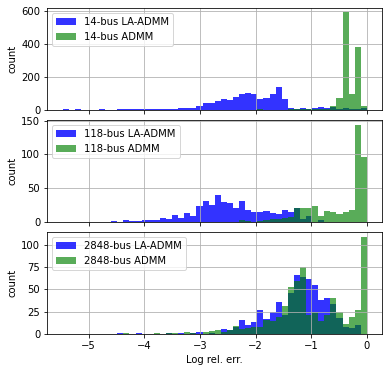}}
\caption{Histograms of $\log_{10}$ relative error in the objective cost at the end of ADMM and LA-ADMM iterations over all test cases.}
\label{fig-cost-hist}
\end{figure}

\subsection{Numerical experiments}
\subsubsection{Test systems}
The three test systems considered were the IEEE 14-bus, IEEE 118-bus, and RTE 2848-bus systems.  Data for all three systems were obtained via the Power Grid Lib (\emph{pglib}) repository \cite{pglib}.  Network properties were obtained directly from the repository, while load data was used to seed the training, as described in the following subsection. 
The partitioning for the 118-bus system was chosen in accordance with the decomposition in \cite{BakerMPC}, which was designed to achieve a reasonable rate of convergence for Lagrangian-based decomposition methods. The RTE 2848-bus network was partitioned into $18$ partitions using spectral factorization of the network graph Laplacian~\cite{hespanha2004efficient}. Network properties for the partitioned systems are summarized in Table \ref{tab-system-data}.

\subsubsection{Simulation of load}
Load for all three systems was sampled such that total load remained less than total generation capacity, but was seeded according to the system \emph{pglib} load data.  In particular, we first defined a characteristic  load for each bus, $\overline{d}_i$, equal to the reported \emph{pglib} value.  We then used a global scaling factor, $\chi$, and bus-level scaling factor, $\xi_i$, to generate the load scenario for each bus via
\begin{align}
\tilde{d}_i(\chi, \xi_i) &= \chi (1 + \xi_i) \overline{d}_i, \quad \chi, \xi_i \in U(0,1), \quad\forall{i}.
\label{eqn-bus-load-sampling}
\end{align}
The random variables $\xi_i\in U(0,1)$ were sampled first, then $\chi$ was sampled uniformly from an interval that was small enough to ensure that the total load did not exceed generation.

\begin{figure}[t!]
\centerline{\includegraphics[width=0.95\linewidth]{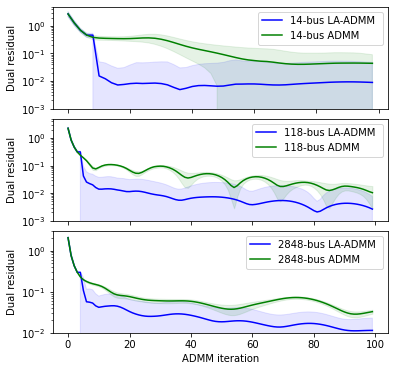}}
\caption{Residual error as function of ADMM iteration averaged over all test cases.  Solid lines indicate the mean while the shaded areas indicate $\pm 1$ standard deviation from the mean.}
\label{fig-aggregate-dual}
\end{figure}

\subsubsection{ADMM configuration}
For all experiments we initialized the ADMM values $\boldsymbol{\lambda}_{su}^{(0)}$ and $\boldsymbol{\theta}_{su}^{(0)}$ to zero.  The network partitions were used to define public buses and corresponding primal and dual variables whose cardinality is given in Table \ref{tab-system-data}.  While a reasonably effective value of $\rho=100$ was chosen for the 118- and 2848-bus systems, we intentionally chose a value of $\rho=1$ for the 14-bus system that led to slow ADMM convergence.  This choice enabled us to demonstrate the efficacy of LA-ADMM even for slowly converging iterations.

\subsubsection{Training and Evaluation}
To generate RNN training data, we sampled instances of system load and solved the resulting DC-OPF problem to obtain converged values of the ADMM objective function as well as primal and dual variables.  While we envision optimal values coming from converged ADMM iterations in any real-world scenario, it is also possible to accelerate training on simulated systems by running ADMM for a small number of iterations to gather what will be the input for the RNN and then using a centralized solution to extract prediction targets (i.e. converged Lagrange multipliers and consensus variables).  We used the latter approach, running ADMM for no more than 10 iterations to obtain training inputs and using the centralized solution for  prediction targets.
Training set sizes were 1000, 4000, and 40,000 for the 14-, 118-, and 2848-bus systems, respectively.

\begin{figure}[t!]
\centerline{\includegraphics[width=0.95\linewidth]{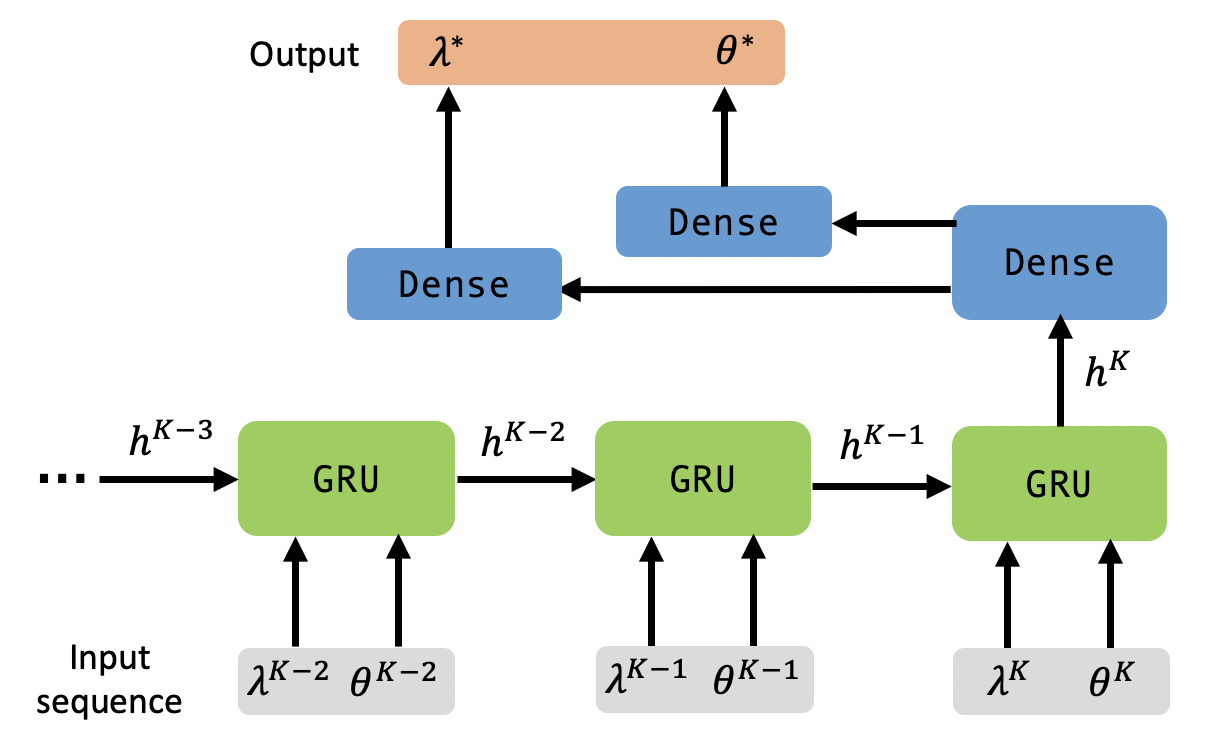}}
\caption{RNN architecture consists of a single Gated Recurrent Unit layer operating on concatenated inputs, followed by dense prediction layers with two-headed linear outputs.}
\label{fig-rnn}
\end{figure}

Best practices for model selection include the use of cross validation, hyperparameter grid search, and regularization, to name a few. After extensive experimentation along these training dimensions, we adopted an approach that favored simplicity over accuracy with respect to a held-out test set.  In particular, we identified a global set of hyperparameters and applied them to the RNN trained for each test system.  The network architecture, illustrated in Fig. \ref{fig-rnn}, consisted of a two-headed input layer corresponding to samples of the form $(\boldsymbol{\lambda}_{su}^{(k),\ell},\boldsymbol{\theta}_{su}^{(k),\ell})_{\ell=1}^L$ for $k=1,...,K$  initial ADMM steps from $\ell=1,\dots,L$ trials.  Inputs were concatenated, then passed through a GRU layer with 128 hidden units, a dense layer of 64 units, and a final two-headed output layer corresponding to targets $(\boldsymbol{\lambda}_{su}^{*,\ell},\boldsymbol{\theta}_{su}^{*,\ell})_{\ell=1}^L$.  All layers except the linear output layer used ReLU activations, and an $L_2$-regularization penalty with coefficient $10^{-4}$ was applied to all weights.   We observed that input sequences of length $K=4$ led to a good balance of data efficiency and prediction accuracy and used this value for all experiments.   Backpropagation was performed with respect to mean squared error on a held-out validation set using the stochastic gradient descent via the Adam optimizer with a learning rate of $10^{-3}$ and stopped after 5 iterations of stagnation in the validation loss or after 50 epochs, whichever came first.

To evaluate performance, new realizations of load were generated using the procedure summarized by \eqref{eqn-bus-load-sampling} and ADMM was run twice per sample: once, uninterrupted, for $100$ iterations to provide a pure ADMM baseline, and a second time with RNN predictions injected at step $k=4$, after which the iteration continues using the standard updates (\ref{eq-admm-xpdate})-(\ref{eq-admm-lambda-update}). 

\begin{figure}[t!]
\centerline{\includegraphics[width=0.95\linewidth]{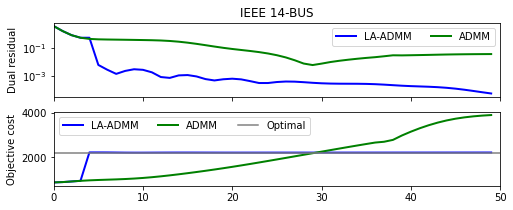}}
\caption{Example of LA-ADMM acceleration for 14-bus test case with 2 partitions and 24 consensus variables.}
\label{fig-14-bus}
\end{figure}

\begin{figure}[t!]
\centerline{\includegraphics[width=0.95\linewidth]{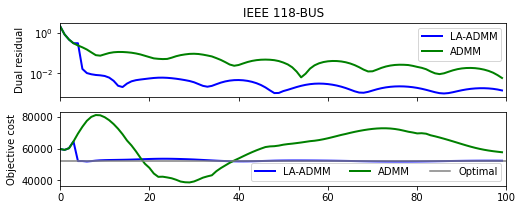}}
\caption{Example of LA-ADMM acceleration for 118-bus test case with 4 partitions and 88 consensus variables.}
\label{fig-118-bus}
\end{figure}

\subsection{Results}
The histograms in Fig. \ref{fig-cost-hist} summarize  convergence of ADMM and LA-ADMM with respect to the known optimal value of the objective function by binning the relative error in this quantity at the end of the iteration.  We see that for the 14- and 118- bus test cases, LA-ADMM almost universally improves convergence by around two orders of magnitude over standard ADMM.  The results for the 2848-bus test case are less striking, although even here we observe acceleration for the slowest converging samples.  Fig. \ref{fig-aggregate-dual} illustrates the effect of injecting predicted values on the mean residual error over all test samples.  Here we observe 1-2 orders of magnitude reduction in the residual on average but that this comes at the cost of higher variance.  Taken together, however, Fig. \ref{fig-cost-hist}-\ref{fig-aggregate-dual} suggest that this is a compelling, data-driven approach to accelerating ADMM.

\begin{figure}[t!]
\centerline{\includegraphics[width=0.95\linewidth]{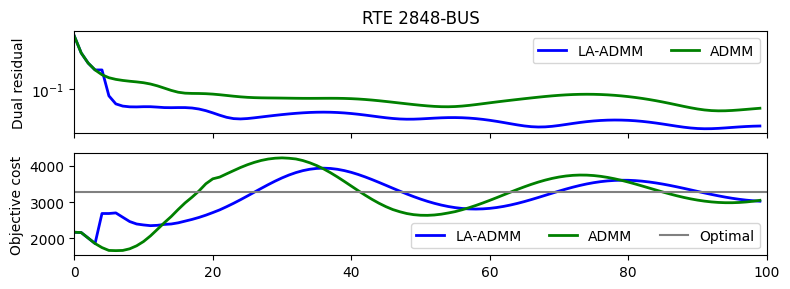}}
\caption{Example of LA-ADMM acceleration for 2848-bus test case with 18 partitions and 1024 consensus variables.}
\label{fig-2848-bus}
\end{figure}

Fig. \ref{fig-14-bus}-\ref{fig-2848-bus} provide a more detailed view of the convergence of ADMM and LA-ADMM for individual test samples.  Fig. \ref{fig-14-bus}-\ref{fig-118-bus} suggest that the acceleration seen in Fig. \ref{fig-cost-hist}-\ref{fig-aggregate-dual} arise from LA-ADMM's ability to push the iteration very close to the optimal solution using predicted values from the RNN.  Analogously, Fig. \ref{fig-2848-bus} shows that LA-ADMM predictions are not as effective in this case.  We speculate that the higher complexity of the modeled system and high dimensionality of the learning problem contribute to LA-ADMM being less effective in this case, but that the model could be improved by increasing the size of the training set and RNN model capacity.  We hope to explore these ideas in future work.

\section{Discussion, Conclusion and Future Work}\label{sec-discussion}
While not necessary apparent today, the use of distributed optimization techniques such as ADMM are likely to be critical in  future power systems operation.  However, such decomposition comes with the drawback of the need to ``iterate to self-consistency".  Our approach avoids this apparent jam by ``pre-optimizing"; while there really is no free lunch, so we must do the work somewhere,  by using computational cycles in an off-line learning mode we allow the online optimization for general inputs to be drastically sped up. 

Generalizability: One may wonder not only whether our method generalizes well to all possible load vectors $\mathbf{d}$ (our numerical tests suggest that it does) but in what other ways it can generalize.  For example, one can imagine both the network structure and generation/flow constraints changing (for the former, if at a day ahead level, say, certain units were not planned to be operating; for the latter, due to, say weather or maintenance). Machine learning models do not readily generalize outside of the domain they were trained on, so these other ways the problem could change would require retraining, but that is not an insurmountable problem.  We will also need to consider whether and how this method can be transferred to the true, nonlinear optimal power flow formulation.

Learning to Optimize:  The specific technique we have invented for this work is part of a growing realization that fertile ground for application of machine learning is not to \emph{replace} optimization but to \emph{augment} it.  In the present case we are using learning to accelerate online optimization of the same sized problem we learn from.  In other cases this idea has been used to learn "heuristics" that allow for efficient solution of much \emph{larger} problems than they were trained on \cite{khalil2017learning} (this is not out of the question here, but it is beyond the scope of the present paper).  

In any case, a benefit of such on approach that is hard to overemphasize is that the learned model is used \emph{within} a generally convergent optimization scheme.  With or without the learning, convergence provides a global, consistent, reliable measure of the quality of the solution.  This helps alleviate a frequent and justified complaint that basing decisions on data-driven machine learning models is not appropriate in safety-critical systems.

\bibliographystyle{IEEEtran}
\bibliography{main}
\end{document}